\documentclass[10pt]{article}
\usepackage{amssymb}
\usepackage{amscd}
\input amssym.def
\input amssym  
\newtheorem{theorem}{Theorem}[section]
\newtheorem{corollary}{Corollary}[section]

\newtheorem{remark}{Remark}[section]
\newtheorem{definition}{Definition}[section]
\newtheorem{proposition}{Proposition}[section]

\def\:{\mbox{\tiny ${\bullet\atop\bullet}$}}

\newcommand{\bea}{\begin{eqnarray}}

\newcommand{\eea}{\end{eqnarray}}
\newcommand{\nn}{\nonumber \\}
\newcommand{\be}{\begin {equation}}

\newcommand{\ee}{\end{equation}}

\usepackage{amssymb}
\title{{The Rogers--Ramanujan recursion and intertwining operators}}
\date{}
\author{S. Capparelli\footnote{S.C. gratefully acknowledges the
partial support of MIUR (Ministero dell'Istruzione, dell'Universit\`a
e della Ricerca).},
J. Lepowsky\footnote{J.L. and A.M. gratefully acknowledge partial
support from NSF grant DMS-0070800.}
and A. Milas}
\begin{document}
\bibliographystyle{alpha}
\maketitle
\begin{abstract}
\noindent We use vertex operator algebras and intertwining operators
to study certain substructures of standard $A_1^{(1)}$--modules,
allowing us to conceptually obtain the classical Rogers-Ramanujan
recursion.  As a consequence we recover Feigin-Stoyanovsky's character
formulas for the principal subspaces of the level 1 standard
$A_{1}^{(1)}$-modules.
\end{abstract}

\renewcommand{\theequation}{\thesection.\arabic{equation}}
\setcounter{equation}{0}
\section{Introduction}
\noindent Vertex operator constructions of affine Lie algebras can be,
and indeed have been, used to prove (or conjecture) many nontrivial
combinatorial identities.  Among them, the most celebrated are the
classical Rogers-Ramanujan identities: 
\bea \label{rr1}
\prod_{i \geq 0} \frac{1}{(1-q^{5i+1}) (1-q^{5i+4})}&=& \sum_{n\geq
0}\frac{q^{n^2}}{(q)_n},\\ 
\label{rr2}
\prod_{i \geq 0} \frac{1}{(1-q^{5i+2})
(1-q^{5i+3})}&=& \sum_{n\geq 0}\frac{q^{n^2+n}}{(q)_n},
\eea
where
\bea
(q)_n &=& (1-q)(1-q^2)\cdots(1-q^n).
\eea
These identities can be expressed in terms of numbers of
partitions.  For example, the first identity states that the number of
partitions of a nonnegative integer with parts of the form $5i+1$ and
$5i+4$ is equal to the number of partitions such that the difference
between consecutive parts is at least 2; this is the so-called {\em
difference-two condition}.  For the classical history of
these identities see \cite{A}.

A proof of these identities was obtained in \cite{LW2}-\cite{LW3} by
means of twisted vertex operators and ``$Z$-operators'' for level-3
standard $A_1^{(1)}$-modules.  In fact, the left-hand side of
(\ref{rr1}) (the product side) was already known to be the
principally specialized character (associated with the principal
gradation) of the vacuum subspace, with respect to a certain
Heisenberg algebra, of a certain standard $A_1^{(1)}$-module
(\cite{LM}, \cite{LW1}), from the Weyl-Kac character formula
\cite{K1}.  In \cite{LW2}-\cite{LW3}, $Z$-operators $Z(j)$ ($j \in
{\mathbb Z}$) were introduced and certain relations were found among
them, giving a basis of the vacuum subspace consisting of monomials in
the $Z(j)$'s applied to the highest weight vector; at the same time,
this construction extended the twisted vertex operator construction
\cite{LW1} of the basic (level 1) standard $A_1^{(1)}$-modules to all
higher levels.  This basis is identified with partitions by virtue of
the indices $j$, and indeed, consecutive indices $j$ satisfy the
difference-two condition.  In the representation-theoretic study of
Rogers-Ramanujan-type identities, a basis of this sort, of a space of
this sort, has come to be called a ``combinatorial basis'' because
(in this case) the difference-two condition is visible from the basis.
Such a basis has also come to be called a ``fermionic basis,'' since
the difference-two condition is an analogue of the difference-one
condition, which corresponds to the classical fermionic statistics of
the Pauli exclusion principle; the difference-two condition is in fact
much more subtle than the difference-one condition, reflecting the
fact that the Rogers-Ramanujan identities (and their generalizations)
are very subtle.
The statement that the ``character'' of (in this original case) the
vacuum space with respect to the Heisenberg algebra of the relevant
level-3 standard $A_1^{(1)}$-module equals the right-hand side (the
sum side) of (\ref{rr1}) has correspondingly come to be called a
``fermionic character formula.''

Analogues of these results were established in \cite{LP} for untwisted,
rather than twisted, vertex operators and $Z$-operators, this time
yielding untwisted $Z$-operator constructions, and in fact,
combinatorial (fermionic) bases, of the standard $A_1^{(1)}$-modules.
In turn, this gave what came to be called ``fermionic character formulas''
analogous to, but different
{}from, the Rogers-Ramanujan identities.  This work extended the
untwisted vertex operator construction of the basic modules for
$A_1^{(1)}$ (\cite{FK}, \cite{S}) to all the higher level modules.

More recently, Feigin and Stoyanovsky \cite{FS1} found
another circle of ideas that led to the Rogers-Ramanujan identities
themselves, using the viewpoint of Lie algebras and loop groups rather
than that of vertex operators.  They observed that
what they called the ``principal subspace''
$W(\Lambda_0)$ (defined in Section 3 below) of the basic $A_1^{(1)}$-module
$L(\Lambda_0)$ has a combinatorial basis that satisfies the
difference-two condition, giving the sum side of (\ref{rr1}).  (This
use of the term ``principal'' is unrelated to the principal gradation
mentioned above.)  However, for principal subspaces there is no
straightforward analogue of the Weyl-Kac character formula and, in
order to get the product side in (\ref{rr1}), Feigin and Stoyanovsky
announced a result that, as a consequence, computes, for all levels,
the characters of the principal subspaces via an analogue of
the Atiyah-Bott fixed point theorem for an infinite-dimensional flag
manifold.  The combinatorial bases obtained in \cite{FS1}-\cite{FS2}
were also implicitly obtained, for all levels for $A_1^{(1)}$,
by Meurman-Primc \cite{MP}, using vertex operator algebra theory.

In \cite{G}, Georgiev extended the character formulas obtained in
\cite{FS1} to a family of standard $A_n^{(1)}$-modules.  This was done
by the explicit construction of combinatorial spanning sets and bases
for the principal subspaces.  In order to prove the linear
independence of the spanning sets, Georgiev used the theory of vertex
operator algebras, including certain intertwining operators.

The present paper is motivated by the fact that many classical, and
also recent, proofs and treatments of the Rogers-Ramanujan identities
and generalizations are based not on the difference-two condition
itself, but rather on a certain classical recursion formula (see \cite{A}):
Noting that the sum sides of the two Rogers-Ramanujan identities can
be obtained as specializations of the generating function
\be
F(x,q)=\sum_{n\geq 0}\frac{x^n q^{n^2}}{(q)_n}
\ee
for $x=1$ and $x=q$, we observe that this generating function satisfies
\be \label{rrrec}
F(x,q)=F(xq,q)+xqF(xq^2,q).
\ee 
We will call this formula the {\em Rogers-Ramanujan recursion} \cite{RR};
see \cite{A} for the role of this recursion in Rogers's and
Ramanujan's work.  
So far, this recursion has not appeared in a
fundamental way in the
vertex-operator studies of the Rogers-Ramanujan identities and related
identities, and we wanted to try to understand it conceptually from
the vertex-operator point of view.

Here is our main result:  We show,
without constructing a basis, that the character of $W(\Lambda_0)$ satisfies
(\ref{rrrec}).
We do this setting up an exact sequence  
$$0 \longrightarrow W(\Lambda_1)
\stackrel{e^{\alpha/2}}{\longrightarrow} W(\Lambda_0)
\stackrel{o(e^{\alpha/2})}{\longrightarrow} W(\Lambda_1)
\longrightarrow 0.$$
Here $W(\Lambda_1)$ is the principal subspace of the second of the
level 1 standard $A_1^{(1)}$-modules (see Section 3 below), and
$e^{\alpha/2}$ and $o(e^{\alpha/2})$ are
the constant {\em factor} and the constant {\em term}, respectively,
of a certain
intertwining operator (cf. \cite{FHL})
associated to $L(\Lambda_0)$ viewed as a vertex operator algebra
(cf. \cite{FLM}).  As a consequence, we recover Feigin-Stoyanovsky's
character formulas for $W(\Lambda_0)$ and $W(\Lambda_1)$.
In \cite{G}, Georgiev had used the same intertwining operator to prove
that the appropriate difference-two combinatorial spanning subset of
$W(\Lambda_0)$ is a basis.

The theory of vertex operator algebras and intertwining operators
(cf. \cite{FLM}, \cite{FHL}, \cite{DL}) is a rich subject, and the
approach initiated here can be modified to apply to affine Lie
algebras of higher rank and level.  This will be treated in separate
publications. In particular, in \cite{CLM} we have extended 
several results from this paper in connection 
with recursions of Rogers and Selberg.

\renewcommand{\theequation}{\thesection.\arabic{equation}}
\setcounter{equation}{0}

\section{Vertex operator constructions associated with $A_{1}^{(1)}$}

In this section we recall some basic definitions and constructions
needed in this paper, especially, the vertex operator construction of
the distinguished basic $A_{1}^{(1)}$-module (\cite{FK}, \cite{S}) and
of its associated vertex operator algebra structure (\cite{B},
\cite{FLM}), and the construction \cite{DL} of a distinguished intertwining
operator associated with its irreducible modules.
 
Let ${{\goth g}}=\goth{sl}(2,\mathbb{C})$ be the 3-dimensional complex
simple Lie algebra with a standard basis $\{ h,\ x_{\alpha},\
x_{-\alpha} \}$ such that $[h,x_{\alpha}]=2x_{\alpha}$,
$[h,x_{-\alpha}]=-2x_{\alpha}$ and $[x_{\alpha},x_{-\alpha}]=h$, and
take the Cartan subalgebra ${\goth{h}}=\mathbb{C} h$.  The standard
symmetric invariant nondegenerate bilinear form $\langle x,y
\rangle={\rm tr}(xy)$ for $x,y \in {{\goth g}}$ allows us to identify
${\goth{h}}$ with its dual ${\goth{h}}^*$.  Take $\alpha \in
{\goth{h}}$ to be the root corresponding to the root vector
$x_{\alpha}$, and take this root to be positive (that is, simple);
then $ \langle \alpha,\alpha \rangle=2$ and we have the root space
decomposition
$${\goth{g}}={\goth{n}}_- \oplus {\goth{h}} \oplus {\goth{n}}_+ ,$$
where ${\goth{n}}_\pm=\mathbb{C}x_{\pm \alpha}$.  Note that under our
identifications,
$$h = \alpha.$$

Consider the Lie algebra
$${\goth{sl}(2)}\ \widehat{} \ =\goth{sl}(2,\mathbb{C}) \otimes
\mathbb{C}[t,t^{-1}] \oplus \mathbb{C}c,$$
where $c$ is a nonzero central element and
$$[x \otimes t^m,y \otimes t^n]=[x,y] \otimes t^{m+n}+\langle x,y
\rangle m \delta_{m+n,0}c$$
for $x,y \in {\goth{g}}$, $m,n \in {\mathbb Z}$.  Adjoining the degree
operator $d$ such that $[d,x \otimes t^m]=m x \otimes t^m$ and
$[d,c]=0$ gives us the Lie algebra $\goth{sl}(2)\
\widetilde{}={\goth{sl}(2)}\ \widehat{} \ \oplus \mathbb{C}d$, the
affine Kac-Moody algebra $A_1^{(1)}$ (cf. \cite{K2}).  
(In this paper we will not need to use $d$.)
Consider the
subalgebra
$$\hat{\goth{h}}_{\mathbb Z}=\coprod_{m \in \mathbb{Z} \setminus \{0\}}
\goth{h} \otimes
t^m \oplus \mathbb{C}c,$$
a Heisenberg algebra in the sense that its commutator subalgebra is
equal to its center, which is one-dimensional, and also consider the
subalgebras
\bea \label{hhat}
\hat{\goth{h}}&=&\goth{h} \otimes \mathbb{C}[t,t^{-1}] \oplus
\mathbb{C}c, \\
\label{nhat}
\goth{n}\widehat{}_+ &=&\goth{n}_+ \otimes \mathbb{C}[t,t^{-1}].
\eea
Throughout the rest of the paper we will
write $x(m)$ for the action of $x \otimes t^m$ on an $\goth{sl}(2)
\ \widehat{} \ $-module, and
we will use the same notation for $\hat{\goth{h}}$-modules.

Now we review the untwisted vertex operator construction of the basic
modules obtained in \cite{FK} and \cite{S} (cf. \cite{FLM}).

Let $P=\frac{1}{2}\mathbb{Z}\alpha$ be the
weight lattice and $Q=\mathbb{Z} \alpha$ the root lattice.
Let $\mathbb{C}[P]$ and $\mathbb{C}[Q]$ be the corresponding group algebras,
with bases $\{e^{\mu} \ | \ \mu \in P\}$ and $\{e^{\mu} \ | \ \mu \in Q\}$.

Consider the induced $\hat{\goth{h}}$-module
$$M(1)=U( \hat{\goth{h}} ) \otimes_{U(\goth{h} \otimes
\mathbb{C}[t] \oplus \mathbb{C}c)} \mathbb{C},$$
where $\goth{h} \otimes \mathbb{C}[t]$ acts trivially on the
one-dimensional module $\mathbb{C}$ and $c$ acts as $1$.
Observe that $M(1)$ is isomorphic to a polynomial algebra in the
infinitely many variables $h(-n)$, $n >0$, that the operators $h(-n)$
act on this polynomial algebra as multiplication operators, that the
operators $h(n)$ act as certain derivations, and that $c$ acts as 1.
Define the vector spaces
$$V_P=M(1) \otimes \mathbb{C}[P],$$
$$V_Q=M(1) \otimes \mathbb{C}[Q],$$
$$V_{Q+\alpha/2}=M(1) \otimes e^{\alpha/2} \mathbb{C}[Q],$$
so that $V_P=V_Q \oplus V_{Q+\alpha/2}$. We define an
$\hat{\goth{h}}$--module structure on $V_P$ by  
making $\hat{\goth{h}}_{\mathbb Z}$ act as $\hat{\goth{h}}_{\mathbb Z}
\otimes 1$ on $V_P=M(1) \otimes \mathbb{C}[P]$ and by making
$\goth{h}=\goth{h} \otimes t^0$ act as $1 \otimes \goth{h}$, with
$h(0)$ defined by
\be \label{action}
h(0)e^\lambda=\langle h,\lambda \rangle e^\lambda
\ee
for $\lambda \in P$.

Let $x, x_0, x_1, x_2 \dots$ 
be commuting formal variables. For $\lambda, \mu \in P$
let $$x^{\lambda} \cdot e^{\mu}=x^{\langle \lambda,\mu \rangle}
e^{\mu}.$$
We extend the action of $x^{\lambda}$ to
$\mathbb{C}[P]$ by linearity and to the whole space
$V_P$ by the formula
$$x^\lambda=1 \otimes x^\lambda.$$
Note that $x^{\lambda}$ is an ${\rm End} \ V_P$-valued formal Laurent
series in the formal variable $x^{1/2}$.
Also let
$$e^\lambda=1 \otimes e^\lambda,$$
acting on $V_P$.
For every $\lambda \in P$, set
\be \label{1.1}
Y(1 \otimes e^{\lambda},x)=E^-(-\lambda,x)E^+(-\lambda,x)
e^{\lambda} x^{\lambda},
\ee
an ${\rm End} \ V_P$-valued formal Laurent series in $x^{1/2}$,
where
$$E^-(\lambda,x)={\rm exp} \left( \sum_{n \leq -1}\frac{\lambda(n)}{n} 
x^{-n} \right),$$
$$E^+(\lambda,x)={\rm exp} \left( \sum_{n \geq 1} \frac{\lambda(n)}{n}
x^{-n} \right).$$
More generally, for the generic homogeneous vector
$w \in V_P$ given by
\be \label{genericvector}
w=h(-n_1-1) \cdots h(-n_k-1) \otimes e^{\lambda}
\ee
with $n_1,\dots,n_k \geq 0$, we set
\be \label{1.1a}
Y(w,x)= \: \prod_{i=1}^k 
\left\{\frac{1}{n_i !} \left(\frac{d}{d x} \right)^{n_i} h(x) \right\} 
Y(1 \otimes e^{\lambda},x) \:,
\ee
where $h(x)=\sum_{n \in \mathbb{Z}} h(n)x^{-n-1}$ and $\: \cdot \:$ 
stands for the normal ordering operation, which places the operators 
$\lambda(n)$ with nonnegative $n$ to the right and 
with negative $n$ to the left.  This formula indeed determines a
well-defined linear map from $V_P$ to the vector space of 
${\rm End} \ V_P$-valued formal Laurent series in $x^{1/2}$ 
(cf. \cite{FLM}). 

We define $\Lambda_0, \Lambda_1 \in (\goth{h} \oplus
\mathbb{C}c)^*$ by:
$\langle \Lambda_i,c \rangle =1$, $\langle \Lambda_i, h
\rangle=\delta_{i,1}$, $i=0,1$. 

\begin{theorem} \label{fks} {\em (\cite{FK}, \cite{S})}
The vector space $V_P$ carries an ${\goth{sl}}(2)\ \widehat{} \
$-module structure uniquely determined by the condition
that for $m \in \mathbb{Z}$, the action $x_{\pm \alpha}(m)$ of
$x_{\pm \alpha} \otimes t^m$ is given by
the coeffcient of $x^{-m-1}$ in $Y(1 \otimes e^{\pm \alpha},x)$.
This module has level 1 (that is, $c$ acts as 1).
Moreover, the direct summands $V_Q$ and $V_{Q+\alpha/2}$ of $V_P$ are
the basic (irreducible) ${\goth{sl}}(2)\ \widehat{} \ $-modules
with highest weights $\Lambda_0$ and
$\Lambda_1$ and with highest weight vectors
$1=1 \otimes 1$ and $e^{\alpha/2}=1 \otimes e^{\alpha/2}$, respectively.
\end{theorem}

In the notation of \cite{K2},
$$V_Q \cong L(\Lambda_0), \ \ V_{Q+\alpha/2} \cong L(\Lambda_1).$$

The following is well known and not hard to prove:

\begin{proposition} \label{x20}
Taking $\lambda = \alpha$ in (\ref{1.1}), we have that the square
$Y(1 \otimes e^{\alpha},x)^2$ of the operator $Y(1 \otimes
e^{\alpha},x)$ on $V_P$ is well defined, and
\be
Y(1 \otimes e^{\alpha},x)^2=0.
\ee
\end{proposition}

It is easy to see that for $m \in {\mathbb Z}$ and $\mu \in P$,
\be \label{1.11a}
x_{\alpha}(m) e^{\mu}=e^{\mu}x_{\alpha}(m+ \langle \alpha,\mu
\rangle),
\ee
and in particular,
\be \label{1.11b}
x_{\alpha}(m) e^{\alpha}=e^{\alpha}x_{\alpha}(m+2)
\ee
and
\be \label{1.11c}
x_{\alpha}(m) e^{\alpha/2}=e^{\alpha/2}x_{\alpha}(m+1).
\ee

We will use the notions of vertex operator algebra, and of module
for such a structure, as defined in \cite{FLM} and \cite{FHL}
(cf. \cite{B}).  We recall the definitions:

A {\em vertex operator algebra} is a vector space $V$ equipped, first, with
a {\em vertex operator map}
\be \label{1.2a}
Y(\cdot ,x) : V \longrightarrow ({\rm End} \ V)[[x,x^{-1}]],
\ee
satisfying the {\em truncation condition}:  For $u,v \in V$, 
the formal Laurent series $Y(u,x)v$ is truncated from below.
It is also equipped with a $\mathbb{Z}$-grading
that is truncated from below, and the homogeneous subspaces are
finite-dimensional.  In addition, $V$ has a {\em vacuum vector}
$\bf{1}$ and a {\em conformal vector} $\omega$.
The main axiom is the {\em Jacobi identity}:
\bea \label{1.3}
\lefteqn{x_0^{-1} \delta \left ( \frac {x_1-x_2}{x_0} \right ) 
Y(u, x_1)Y(v, x_2)} \nn
&& -x_0^{-1} \delta \left ( \frac
{x_2-x_1}{-x_0} \right ) Y(v, x_2) 
Y(u, x_1)
\nn
&&= x_2^{-1} \delta \left ( \frac {x_1-x_0}{x_2} \right )Y(Y(u, x_0)v, x_2)
\eea
for $u,v \in V$.  Here 
$$\delta(x)=\sum_{n \in \mathbb{Z}} x^n$$
is the formal delta function, and in the expansions of the three
delta-function expressions, the negative powers of the binomials are
understood to be expanded in nonegative powers of the second variable.
In addition, $Y({\bf{1}},x)$ is the identity operator on $V$,
and the {\em creation property} holds: For $v \in V$, $Y(v,x){\bf{1}}$
involves only nonnegative powers of $x$ and its constant term is $v$.
Also,
$$Y(\omega,x) = \sum_{n \in \mathbb{Z}} L(n)x^{-n-2},$$
where the operators $L(n)$ on $V$ satisfy the standard bracket
relations for the Virasoro algebra (including a central charge);
the $\mathbb{Z}$-grading on $V$ coincides with the eigenspace
decomposition of the operator $L(0)$; and the operator $L(-1)$
satisfies the {\it $L(-1)$-derivative property}:
$$Y(L(-1) v, x)=\frac{d}{d x} Y(v, x)$$
for $v \in V$.

A {\em module}
for the vertex operator algebra $V$ is a vector space $W$ 
equipped with a {\em vertex operator map}
\be
Y(\cdot ,x) : V \longrightarrow ({\rm End} \ W)[[x,x^{-1}]],
\ee
such that all the axioms in the definition of the notion of vertex
operator algebra that make sense hold, except that the grading on $W$
is allowed to be a $\mathbb{Q}$-grading rather than a
$\mathbb{Z}$-grading.

We have already encountered some special vertex operators in
(\ref{1.1}) and (\ref{1.1a}), in the sense of the following
result (\cite{B}, \cite{FLM}):

\begin{theorem} \label{voa}
Formulas (\ref{1.1}) and (\ref{1.1a}) give
a vertex operator algebra structure on $V_Q$, and
$V_{Q+\alpha/2}$ is a $V_Q$-module.  The vertex operator algebra
$V_Q$ is simple and the module $V_{Q+\alpha/2}$ is irreducible.
\end{theorem}

Even though the vertex operators (\ref{1.1}) and (\ref{1.1a}) are
defined for every $w \in V_P$ (recall the comment after (\ref{1.1a})),
the vertex operator algebra structure
$V_Q$ cannot be extended to such a structure on $V_P$; instead, $V_P$ has 
the structure of an {\em abelian intertwining algebra} in the
sense of \cite{DL}.
However, following \cite{DL}, for $w \in V_{Q+\alpha/2}$
we modify the vertex operator map $Y$ by:
\be \label{modifiedY}
{\mathcal Y}(w,x)=Y(w,x)e^{i \pi \alpha/2};
\ee
then a version of the Jacobi identity (\ref{1.3}) still holds \cite{DL}:
\bea \label{1.4}
\lefteqn{x_0^{-1} \delta \left ( \frac {x_1-x_2}{x_0} \right ) 
Y(u, x_1)\mathcal{Y}(w, x_2)} \nn
&& -x_0^{-1} \delta \left ( \frac
{x_2-x_1}{-x_0} \right ) \mathcal{Y}(w, x_2) 
Y(u, x_1)
\nn
&&= x_2^{-1} \delta \left ( \frac {x_1-x_0}{x_2} \right
)\mathcal{Y}(Y(u, x_0)w, x_2)
\eea
for $u \in V_Q$.
In addition, it is not hard to see that
\be \label{1.5}
\mathcal{Y}(w,x)|_{V_Q} {V_Q} \subset V_{Q+\alpha/2}((x))
\ee
\be  \label{1.6}
\mathcal{Y}(w,x)|_{V_{Q+\alpha/2}} V_{Q+\alpha/2} \subset V_{Q}((x^{1/2})),
\ee
where, for a vector space $W$ and a formal variable $y$, $W((y))$
denotes the space of $W$-valued lower truncated formal Laurent series
in $y$.

Also for a general a vector space $W$, we define $W \{ x \}$ to be the
vector space of $W$-valued formal expressions of the form $\sum_{s \in
\mathbb{Q}} w_s x^s$, $w_s \in W$.  We recall from \cite{FHL} the
definition of the notion of intertwining operator for a triple of
modules for a vertex operator algebra:

\begin{definition} \label{intertwining}
{\em Let $W_1,W_2$ and $W_3$ be modules for a vertex
operator algebra $V$. A map
$${\cal Y}( \cdot ,x) \ :  \ W_1 \longrightarrow ({\rm Hom}(W_2,W_3) )
\{x \},$$
is an {\em intertwining operator} of {\em type} ${W_3 \choose W_1 \
W_2}$ if it satisfies all the defining properties of a module action
that make sense, namely,
the truncation property (that for $w_i
\in W_i$, $i=1,2$, $\mathcal{Y}(w_1,x)w_2$ is truncated from below),
the Jacobi identity
\bea
\lefteqn{x_0^{-1} \delta \left ( \frac {x_1-x_2}{x_0} \right ) 
Y(u, x_1) \mathcal{Y}(w_1, x_2)w_2} \nn
&& -x_0^{-1} \delta \left ( \frac
{x_2-x_1}{-x_0} \right ) \mathcal{Y}(w_1, x_2) 
Y(u, x_1)w_2
\nn
&&= x_2^{-1} \delta \left ( \frac {x_1-x_0}{x_2} \right
)\mathcal{Y}(Y(u, x_0)w_1, x_2)w_2
\eea
for $u\in V$, $w_1 \in W_1$ and $w_2 \in W_2$,
and the $L(-1)$-derivative property: for $w_1 \in W_1$, 
$$\mathcal{Y}(L(-1) w_1, x)=\frac{d}{d x}
\mathcal{Y}(w_1, x).$$
We will denote the vector space of all intertwining operators of type 
${W_3 \choose W_1 \  W_2}$
by $I \ { W_3 \choose W_1 \  W_2}$.}
\end{definition}

The operator (\ref{modifiedY}) in fact gives intertwining operators:

\begin{theorem} \label{intwop} {\em (\cite{DL}, \cite{FHL})}
The operators (\ref{1.5}) and (\ref{1.6}) are 
intertwining operators of the types
$\displaystyle{{L(\Lambda_1) \choose L(\Lambda_1) \ L(\Lambda_0) }}$ and
$\displaystyle{{L(\Lambda_0) \choose L(\Lambda_1) \ L(\Lambda_1) }}$,
respectively.  Moreover, the spaces
$I \displaystyle{{L(\Lambda_1) \choose L(\Lambda_1) \  L(\Lambda_0)}}$ and
$I \displaystyle{{L(\Lambda_0) \choose L(\Lambda_1) \  L(\Lambda_1)
}}$ are one-dimensional.
\end{theorem}

\renewcommand{\theequation}{\thesection.\arabic{equation}}
\setcounter{equation}{0}

\section{The principal subspaces}

We recall the definition of the {\em principal subspaces}
$W(\Lambda_0) \subset L(\Lambda_0)$ and 
$W(\Lambda_1) \subset L(\Lambda_1)$ from \cite{FS1}:
$$W(\Lambda_0)= U({\goth{n}}\widehat{}_+) v_{\Lambda_0},$$  
$$W(\Lambda_1)= U({\goth{n}}\widehat{}_+) v_{\Lambda_1},$$ 
where 
$v_{\Lambda_0}$ and $v_{\Lambda_1}$ are highest weight vectors for
$L(\Lambda_0)$ and $L(\Lambda_1)$, respectively.

Clearly, $W(\Lambda_j)$, $j=0,1$,
is spanned by
\be \label{combspan}
x_{\alpha}(-m_1) \cdots x_{\alpha}(-m_k)v_{\Lambda_j},
\ee
where $m_1, \dots ,m_k > 0$.  
{}From now on we will identify 
$L(\Lambda_0)$ with $V_{Q}$ and $L(\Lambda_1)$ with $V_{Q+\alpha/2}$ 
and we will take
$$v_{\Lambda_0}=1, \ \ v_{\Lambda_1}=e^{\alpha/2}v_{\Lambda_0}=e^{\alpha/2}.$$

The action of $L(0)$ on the space $V_P$ gives us the grading of $V_P$
by {\em weights}, determined by:
\be
{\rm wt} \ e^{\lambda} = \frac {1}{2} \langle \lambda, \lambda \rangle
\ee
for $\lambda \in P$ and by the condition that for $n \in \mathbb{Z}$,
the weight of the operator $h(-n)$ on $V_P$ is $n$; this determines
the weight of the generic vector (\ref{genericvector}).  This grading
by weights is a $\mathbb{Q}$-grading, or more specifically, a $\frac
{1}{4} \mathbb{Z}$-grading.  The space $V_P$ also has a second,
compatible, grading, by {\em charge}, given by the eigenvalues of
the operator $\frac{1}{2}\alpha(0)=\frac{1}{2}h(0)$
(recall (\ref{action})); this a 
$\frac {1}{2} \mathbb{Z}$-grading.  We shall consider these gradings
restricted to the principal subspaces $W(\Lambda_0)$ and
$W(\Lambda_1)$.  The formulas
$$x^{\alpha/2} x_{\alpha}(-m_1) \cdots x_{\alpha}(-m_k)v_{\Lambda_0}=
x^k x_{\alpha}(-m_1) \cdots x_{\alpha}(-m_k)v_{\Lambda_0},$$
$$x^{\alpha/2} x_{\alpha}(-m_1) \cdots x_{\alpha}(-m_k)v_{\Lambda_1}=
x^{k+1/2} x_{\alpha}(-m_1) \cdots x_{\alpha}(-m_k)v_{\Lambda_1}$$
show that the charges of the indicated elements of
$W(\Lambda_0)$ and $W(\Lambda_1)$ are 
$k$ and $k+1/2$, respectively.

For a subspace $M$ of $V_P$ homogeneous with respect to the double
grading, we consider the corresponding graded dimension of $M$---the
generating function of the dimensions of its homogeneous subspaces,
where we use the formal variables $x$ and $q$:
\be
{\rm dim}_*(M;x,q) = \mbox{tr}|_{M}x^{\alpha/2}q^{L(0)}.
\ee
We shall focus on the graded dimensions of the principal subspaces 
$W(\Lambda_0)$ and $W(\Lambda_1)$, which we shall write as:
\bea
\label{9.1}
&&\chi_{0}(x,q) = {\rm dim}_*(W(\Lambda_0);x,q)
= \mbox{tr}|_{W(\Lambda_0)}x^{\alpha/2}q^{L(0)}, \\ 
&&\chi_{1}(x,q) = {\rm dim}_*(W(\Lambda_1);x,q)
= \mbox{tr}|_{W(\Lambda_1)}x^{\alpha/2}q^{L(0)}.
\eea

One way to describe the spaces $W(\Lambda_0)$ and $W(\Lambda_1)$ is
to use the vertex operators
\be \label{9.4}
Y(e^{\alpha},x_i)=E^-(-\alpha,x_i)E^+(-\alpha,x_i)e^{\alpha}x_i^{\alpha},
\ee
(we identify  $e^{\alpha}$ with $1 \otimes e^{\alpha}$)
or the slightly modified operators
\be \label{9.5}
X(e^{\alpha},x_i)=E^-(-\alpha,x_i)E^+(-\alpha,x_i)e^{\alpha}x_i^{\alpha+1}.
\ee
The operators $Y(e^{\alpha},x_i)$ and $X(e^{\alpha},x_i)$ all commute
with one another, and
\bea \label{9.6}
&& \prod_{i=1}^n X(e^{\alpha},x_i)v_{\Lambda_0}= \nn
&& =\prod_{i < j} (x_i-x_j)^2 \prod_{i=1}^n x_i \: \prod_{i=1}^n
E^-(-\alpha,x_i) 
E^+(-\alpha,x_i) \: e^{n \alpha}
\nn
&& = \prod_{i<j}(x_i-x_j)^2 \prod_{i=1}^n x_i
\prod_{i=1}^n E^-(-\alpha,x_i) e^{n\alpha}.
\eea
Similarly, for $W(\Lambda_1)$,
\bea \label{9.7}
&& \prod_{i=1}^n X(e^{\alpha},x_i)v_{\Lambda_1}=\prod_{i=1}^n
X(e^{\alpha},x_i)e^{\alpha/2} \nn
&& = \prod_{i<j}(x_i-x_j)^2 \prod_{i=1}^n x_i^2
\prod_{i=1}^n E^-(-\alpha,x_i)
e^{(n+1/2)\alpha}.
\eea

Another way to describe the principal subspaces is as certain
quotients of the associative algebra
$U({\goth{n}}\widehat{}_+ )$, which
is in fact commutative.
This description will be used from now on,
and we will not directly use the formulas (\ref{9.6})
and (\ref{9.7}), except to compute the action of components of the
operator (\ref{9.4}) (or (\ref{9.5})) on particular vectors.

Consider a set of commuting independent formal variables
$y_{-1}, y_{-2}, \ldots $ and the corresponding polynomial algebra 
$$\mathcal A= \mathbb{C}[y_{-1},y_{-2}, \ldots ].$$
Consider the algebra map
$$\mathcal{A} \longrightarrow {\rm End} \ V_P$$
$$y_{-j} \mapsto x_{\alpha}(-j)$$
($j>0$); this map is well defined because the operators 
$x_{\alpha}(-j)$ commute. We shall use this
correspondence to describe $W(\Lambda_0)$ and
$W(\Lambda_1)$.
Define the linear surjection
\bea \label{9.11}
f_{\Lambda_0}: \mathcal A & \rightarrow & W(\Lambda_0) \nn
p(y_{-1},y_{-2},\dots)
& \mapsto & p(x_{\alpha}(-1),x_{\alpha}(-2),\ldots) \cdot
v_{\Lambda_0}
\eea
and set $$\mathcal A_{\Lambda_0}= {\rm Ker}\ f_{\Lambda_0},$$ an ideal in $\mathcal{A}$.
Then we have
\be \label{9.12}
\mathcal A/\mathcal A_{\Lambda_0}
\stackrel{\sim}{\longrightarrow} W(\Lambda_0). 
\ee

Similarly, for $W(\Lambda_1)$, we have
\bea \label{9.14}
f_{\Lambda_1}: \mathcal A & \rightarrow & W(\Lambda_1)\nn
p & \mapsto & p(x_{\alpha}(-1),x_{\alpha}(-2),\ldots) \cdot
v_{\Lambda_1}
\eea
and we set
$$\mathcal A_{\Lambda_1}= {\rm Ker} \ f_{\Lambda_1}.$$
Then 
$$\mathcal{A}/ \mathcal A_{\Lambda_1} \stackrel{\sim}{\longrightarrow} W(\Lambda_1).$$
Notice that 
$$
y_{-1}\in \mathcal A_{\Lambda_1}
$$ 
since 
\be \label{9.15}
x_{\alpha}(-1)\cdot e^{\alpha/2}v_{\Lambda_0}=0,
\ee
which follows from (\ref{9.7}) for $n=1$.
Since  $\mathcal A_{\Lambda_1}$ is an ideal, it includes
the ideal generated by $y_{-1}$:
$$(y_{-1}) = \mathcal A y_{-1}
\subset \mathcal A_{\Lambda_1}.$$ 
This ideal
$(y_{-1}) $ has a natural complement in $\mathcal{A}$, namely,
\be \label{9.16}
\mathcal A'= \mathbb C[y_{-2},y_{-3},\ldots ];
\ee
that is,
\be \label{9.17}
\mathcal A = (y_{-1}) \oplus \mathcal A'. 
\ee
Thus we have a natural linear surjection
\bea \label{9.18}
f'_{\Lambda_1}=f_{\Lambda_1}|_{\mathcal A'}: \mathcal A' & \rightarrow &
W(\Lambda_1) \nn
p(y_{-2},y_{-3},\dots) & \mapsto & p(x_{\alpha}(-2),x_{\alpha}(-3),\dots)
\cdot v_{\Lambda_1}.
\eea
Define
$$
\mathcal B_{\Lambda_1}= {\rm Ker} f'_{\Lambda_1},
$$
an ideal in $\mathcal{A}'$. Then  
\be \label{9.20}
\mathcal A'/\mathcal B_{\Lambda_1} \stackrel{\sim}{\longrightarrow} W(\Lambda_1).
\ee
We also have the natural surjection (\ref{9.14}),  and its kernel
$\mathcal{A}_{\Lambda_1}$ can be decomposed
as
\be \label{9.22}
\mathcal A_{\Lambda_1}=(y_{-1})\oplus \mathcal B_{\Lambda_1}\subset \mathcal A.
\ee
Thus the map $f_{\Lambda_1}$ induces natural isomorphisms
\be \label{9.23}
 \mathcal A/\mathcal A_{\Lambda_1}= ((y_{-1})\oplus \mathcal
A')/ ( (y_{-1}) \oplus \mathcal B_{\Lambda_1})
\stackrel{\sim}{\longrightarrow}
\mathcal A'/\mathcal B_{\Lambda_1} 
\stackrel{\sim}{\longrightarrow}  W(\Lambda_1).
\ee

Consider the linear map 
\be \label{ealpha/2}
e^{\alpha/2} : V_P \rightarrow V_P.
\ee
This is clearly a linear isomorphism, since its inverse
is $e^{-\alpha/2}$. Let us restrict this map
to $W(\Lambda_0)$. Then we have
\be \label{9.24}
e^{\alpha/2} : W(\Lambda_0) \longrightarrow W(\Lambda_1),
\ee
since, from (\ref{1.11c}),
\be \label{9.24a}
e^{\alpha/2}(x_{\alpha}(-i_1)\cdots x_{\alpha}(-i_k)\cdot 1) =
x_{\alpha}(-i_1-1)\cdots x_{\alpha}(-i_k-1 ) \cdot e^{\alpha/2}
\ee
for $i_j \in \mathbb Z$.
The map (\ref{9.24}) is certainly injective, and it is surjective
because
\be \label{9.24b}
U({\goth{n}}\widehat{}_+) \cdot e^{\alpha/2}=
e^{\alpha/2}(U({\goth{n}}\widehat{}_+) \cdot 1).
\ee
Hence (\ref{9.24}) is a linear isomorphism.

We construct a lifting 
$\widehat{e^{\alpha/2}}$ of $e^{\alpha/2}$ as follows:
\bea \label{9.25}
\widehat{e^{\alpha/2}}:
\mathcal A & \stackrel{\sim}{\rightarrow} & \mathcal A' \nn
y_{-j} & \mapsto & y_{-j-1}
\eea
($j\geq 1$); this is an algebra isomorphism. Then $\widehat{e^{\alpha/2}}$
is a lifting in the sense that the following diagram commutes:
$$
\CD
\mathcal A @> \widehat{e^{\alpha/2}}> \sim > \mathcal A' @>\iota 
>>\mathcal A \\
@Vf_{\Lambda_0}VV @Vf'_{\Lambda_1}VV @Vf_{\Lambda_1}VV \\
W(\Lambda_0)\simeq \mathcal A/\mathcal A_{\Lambda_0} @> {e^{\alpha/2}}>\sim >
W(\Lambda_1)\simeq
\mathcal A'/\mathcal B_{\Lambda_1} @> \sim >> \mathcal A/\mathcal A_{\Lambda_1}
\endCD
$$
where $\iota : \ \mathcal{A}' \hookrightarrow \mathcal{A}$ is the algebra
injection. The right half of this diagram clearly commutes, and we now
verify that the left half commutes: For any monomial 
$$
M=y_{-j_1} \cdots y_{-j_n} \in \mathcal A
$$
($j_p\geq 1$), we have
$$
f_{\Lambda_0}(M)=x_{\alpha}(-j_1) \cdots x_{\alpha}(-j_n)\cdot
v_{\Lambda_0} \in W(\Lambda_0), 
$$
$$
e^{\alpha/2}( f_{\Lambda_0}(M))=x_{\alpha}(-j_1-1) \cdots
x_{\alpha}(-j_n-1)\cdot e^{\alpha/2} \in W(\Lambda_1), 
$$
and also,
$$
\widehat{e^{\alpha/2}}(M)= y_{-j_1-1} \cdots y_{-j_n-1} \in
\mathcal{A}',
$$
$$
f'_{\Lambda_1}(\widehat{e^{\alpha/2}}(M))=x_{\alpha}(-j_1-1) \cdots
x_{\alpha}(-j_n-1)\cdot 
e^{\alpha/2} \in W(\Lambda_1). 
$$
Thus the diagram commutes.

It follows that $\widehat{e^{\alpha/2}}$ maps $\mathcal A_{\Lambda_0}$ isomorphically
onto $\mathcal B_{\Lambda_1}$:
\be \label{9.31}
\widehat{e^{\alpha/2}}:\mathcal A_{\Lambda_0} \stackrel{\sim}{\longrightarrow} \mathcal
B_{\Lambda_1}.
\ee
That is, the relations for $W(\Lambda_1)$ are precisely the
index-shifted relations for
$W(\Lambda_0)$.

It is important to notice that the linear isomorphism 
$e^{\alpha/2}:W(\Lambda_0) \rightarrow W(\Lambda_1)$ does not preserve
charge or weight.
In fact, we will prove that 
\be \label{9.36}
\chi_1 (x,q)= x^{1/2}q^{1/4}\chi_0(xq,q),
\ee
that is,
\be \label{9.37}
 \sum _{r,s} \mbox{dim }W(\Lambda_1)_{r,s}x^{r}q^{s}=\sum _{r,s}
\mbox{dim }W(\Lambda_0)_{r,s}x^{r+1/2}q^{r+s+1/4}, 
\ee
where $W(\Lambda_j)_{r,s}$ denotes the vector subspace of
$W(\Lambda_j)$ of charge $r \in \frac{1}{2}\mathbb{Z}$ and weight
$s \in \frac{1}{4}\mathbb{Z}$.  
This is equivalent to
\be \label{9.38}
x^{-1/2}q^{-1/4}\chi_1(x,q)=\chi_0(xq,q),
\ee
that is, to
\bea \label{9.39}
&& \sum _{r,s}
\mbox{dim }W(\Lambda_1)_{r+1/2,s+1/4}x^{r}q^{s}=
\sum _{r,s} \mbox{dim }W(\Lambda_1)_{r,s}x^{r-1/2}q^{s-1/4}=\nn
&& \sum_{r,s} \mbox{dim }W(\Lambda_0)_{r,s}x^{r}q^{r+s}=
\sum _{r,s} \mbox{dim }W(\Lambda_0)_{r,s-r}x^{r}q^{s}.
\eea
To prove (\ref{9.36}) it is enough to show 
$$\mbox{dim }W(\Lambda_0)_{r,s}=\mbox{dim
}W(\Lambda_1)_{r+1/2, r+s+1/4}.$$
Let us restrict $e^{\alpha/2}$ to $W(\Lambda_0)_{r,s}$.  Then
(\ref{9.24a}) and the fact that
$\frac{1}{2} \langle \frac{\alpha}{2},\frac{\alpha}{2} \rangle =\frac{1}{4}$
show that
$$e^{\alpha/2}: W(\Lambda_0)_{r,s} \longrightarrow W(\Lambda_1)_
{r+1/2, r+s+1/4},$$
and similarly, applying $e^{-\alpha/2}$ to (\ref{9.24a}), we have
$$e^{-\alpha/2}: W(\Lambda_1)_{r+1/2, r+s+1/4}\longrightarrow 
W(\Lambda_0)_{r,s}.$$
Thus we have the isomorphism
\be \label{9.40}
e^{\alpha/2}: W(\Lambda_0)_{r,s} \stackrel{\sim}{\longrightarrow}
W(\Lambda_1)_{r+1/2, r+s+1/4}.
\ee
This proves (\ref{9.36}), which in effect describes precisely how
$e^{\alpha/2} : W(\Lambda_0) \rightarrow W(\Lambda_1)$
shifts charge and weight.

\begin{remark}
Note that $\chi_1(1,q)=q^{1/4}\chi_0(q,q)$ gives the principally
specialized character of $W(\Lambda_0)$, while $\chi_0(1,q)$ gives (by
definition) the ordinary graded dimension of $W(\Lambda_0)$ with
respect to the grading by weights.
\end{remark}

\renewcommand{\theequation}{\thesection.\arabic{equation}}
\setcounter{equation}{0}
\section{The main theorem}

We have:

\begin{theorem} \label{mainlemma} 
Take the intertwining operator
$\mathcal{Y}( \cdot,x)  \in I \ {L(\Lambda_1) \choose L(\Lambda_1) \
L(\Lambda_0)}$ as given in (\ref{modifiedY}) and (\ref{1.5}).
The sequence
\be \label{10.1}
0\longrightarrow W(\Lambda_1) \stackrel{ e^{\alpha/2}}{\longrightarrow}
W(\Lambda_0) \stackrel{o(e^{\alpha/2})}{\longrightarrow} W(\Lambda_1)
\longrightarrow 0 
\ee
is exact, where 
$$o(e^{\alpha/2})={\rm Res}_x x^{-1}\mathcal Y(e^{\alpha/2},x),$$
the constant term of the intertwining operator 
$\mathcal Y(e^{\alpha/2},x)$.
\end{theorem}

\begin{remark}
Note that  $e^{\alpha/2}$ is the ``canonical constant {\it factor}'' of
$\mathcal Y(e^{\alpha/2},x)$, where $\mathcal{Y}( \cdot ,x)$ is viewed as
in (\ref{1.6}), and that $o(e^{\alpha/2})$ is the
``canonical constant {\it summand}" of $\mathcal Y(e^{\alpha/2},x)$,
where $\mathcal{Y}( \cdot ,x)$ is viewed as
in (\ref{1.5}).
\end{remark}
{\em Proof of Theorem \ref{mainlemma}:}
Recall that $e^{\alpha/2}: W(\Lambda_0) \rightarrow W(\Lambda_1)$
is a linear isomorphism; note on the
other hand that the map $e^{\alpha/2}$ in (\ref{10.1})
goes in the reverse direction. We shall see below that $e^{\alpha/2}$
takes $W(\Lambda_1)$ into $W(\Lambda_0)$; it is certainly injective.
Let us recall the Jacobi identity (\ref{1.4}).  Taking
$u=e^{\alpha}$ and $w=e^{\alpha/2}$ in (\ref{1.4})
and applying the residue operation
${\rm Res}_{x_0}$, we find that
\bea
&& [\mathcal Y(e^{\alpha/2},x_1),Y(e^\alpha,x_2)]=0,
\eea
so that
\bea \label{10.1a}
&& [\mathcal Y(e^{\alpha/2},x), {\goth{n}}\widehat{}_+]=0.
\eea
In addition,
$$o(e^{\alpha/2})v_{\Lambda_0}=v_{\Lambda_1},$$
and so
\be \label{10.2}
{\rm Res}_x x^{-1} \mathcal Y(e^{\alpha/2},x) U({\goth{n}}\widehat{}_+)v_{\Lambda_0}=
U({\goth{n}} \widehat{}_+)v_{\Lambda_1}.
\ee
Hence $o(e^{\alpha/2})$ takes $W(\Lambda_0)$ to $W(\Lambda_1)$
and is surjective. 

Now we show that $e^{\alpha/2} W(\Lambda_1) \subset W(\Lambda_0)$
and that 
\be \label{10.3}
e^{\alpha/2}W(\Lambda_1)\subset {\rm Ker} \ o(e^{\alpha/2}),
\ee
so that the sequence (\ref{10.1}) is a chain complex:
For $j_1,...,j_n \in \mathbb{Z}$,
\bea \label{10.5}
&& e^{\alpha/2}x_{\alpha}(j_1)\cdots
x_{\alpha}(j_n)e^{\alpha/2}v_{\Lambda_0}=\nn
&& x_{\alpha}(j_1-1)\cdots x_{\alpha}(j_n-1)e^{\alpha}v_{\Lambda_0}= \nn
&& x_{\alpha}(j_1-1)\cdots x_{\alpha}(j_n-1)x_{\alpha}(-1)v_{\Lambda_0},
\eea
which gives that $e^{\alpha/2} W(\Lambda_1) \subset W(\Lambda_0)$. The fact
that  $o(e^{\alpha/2})e^{\alpha/2} W(\Lambda_1)=0$ now
follows from (\ref{10.1a}) and the fact 
$x_{\alpha}(-1)e^{\alpha/2}=0$.
(Incidentally,
$$x_{\alpha}(-1)e^{\alpha/2}v_{\Lambda_0}=e^{\alpha/2}x_{\alpha}(0)
v_{\Lambda_0}=0;$$
that is, the fact that $x_{\alpha}(-1)$ annihilates
$v_{\Lambda_1}$ is equivalent to the fact that $x_{\alpha}(0)$ annihilates
$v_{\Lambda_0}$, by means of $e^{\alpha/2}$.
So the chain complex property follows from {\it either}
$x_{\alpha}(-1)v_{\Lambda_1}=0$ {\it or}
$x_{\alpha}(0)v_{\Lambda_0}=0$.)
Thus the sequence (\ref{10.1}) is a chain complex, and
it remains to prove the exactness.

We shall now characterize the kernel of $o(e^{\alpha/2})\ :
\ W(\Lambda_0) \rightarrow W(\Lambda_1)$.

Let $p$ be a polynomial in $\mathcal A$, so that $f_{\Lambda_0}(p)$ is a
general element of
$W(\Lambda_0)$. Then 
\bea \label{10.6}
&& o(e^{\alpha/2})(f_{\Lambda_0}(p))=0 \iff  o(e^{\alpha/2})
(p(x_{\alpha}(-1),x_{\alpha}(-2),\ldots) v_{\Lambda_0})=0 \nn
&& \iff p(x_{\alpha}(-1),x_{\alpha}(-2),\ldots) v_{\Lambda_1}=0 \\
&& \iff p\in {\rm Ker} \ f_{\Lambda_1}= \mathcal A_{\Lambda_1}=
(y_{-1})\oplus \mathcal
B_{\Lambda_1}
\iff \pi(p)\in \mathcal B_{\Lambda_1}=\widehat {e^{\alpha/2}}
(\mathcal{A}_{\Lambda_0}), \nonumber
\eea
where $$\pi: \mathcal A \rightarrow \mathcal A'$$ is the projection
with respect to
the decomposition (\ref{9.17}); $\pi$ is an algebra surjection.
Thus
\be \label{10.7}
f_{\Lambda_0}(p) \in {\rm Ker} \ o(e^{\alpha/2})
\iff \pi(p)=\widehat {e^{\alpha/2}}(p')
\ee
for some $p'\in \mathcal A_{\Lambda_0}$,
and so 
$$f_{\Lambda_0}(p)\in {\rm Ker}\ o(e^{\alpha/2}) \iff 
(\widehat {e^{\alpha/2}}^{-1} \circ \pi )(p) \in \mathcal A_{\Lambda_0}.$$

Note that if $p\in \mathcal A_{\Lambda_0}$, then $f_{\Lambda_0}(p)=0$, so
that all the above statements
hold in particular for $p$.
Thus
\be \label{10.8}
p\in \mathcal A_{\Lambda_0} \Rightarrow \pi(p) \in \widehat
{e^{\alpha/2}}(\mathcal A_{\Lambda_0}),
\ee
i.e.,
\be \label{10.9}
\pi(\mathcal A_{\Lambda_0}) \subset \widehat {e^{\alpha/2}}(\mathcal
A_{\Lambda_0})\  (=\mathcal B_{\Lambda_1}),
\ee
i.e.,
\be \label{10.10}
(\widehat {e^{\alpha/2}}^{-1} \circ \pi ) (\mathcal A_{\Lambda_0})\subset
\mathcal A_{\Lambda_0}.
\ee

Define
\be \label{10.11}
S= \widehat {e^{\alpha/2}}^{-1} \circ \pi : \mathcal A \rightarrow \mathcal
A.
\ee
This operator is the algebra surjection, with kernel
$(y_{-1})$, defined by
$$
y_{-j} \mapsto y_{-j+1}
$$
for $j\geq 2$ and
$$
y_{-1}\mapsto 0, 
$$
and thus is a ``shift map''. The map
$S$ is 0 on the ideal $(y_{-1})$ and is the isomorphism
\be \label{10.14}
 \widehat {e^{\alpha/2}}^{-1} : \mathcal A' \rightarrow \mathcal
A
\ee
on the complement $\mathcal{A}'$ of $(y_{-1})$ in $\mathcal{A}$.
With this definition we can restate our characterization of the kernel
as follows: For $p \in \mathcal{A}$,
\be \label{10.15}
f_{\Lambda_0}(p)\in {\rm Ker} \ o(e^{\alpha/2})\iff S(p) \in \mathcal
A_{\Lambda_0}
\left( \iff f_{\Lambda_0}(S(p))=0 \right). 
\ee
We have also noted that
\be \label{10.16}
S(\mathcal A_{\Lambda_0})\subset \mathcal A_{\Lambda_0}
\ee
(i.e., if $f_{\Lambda_0}(p)=0$, then $f_{\Lambda_0}(S(p))=0$).

Next we characterize the image of $e^{\alpha/2}: \ W(\Lambda_1) \rightarrow 
W(\Lambda_0)$: For $p \in \mathcal{A}$, when is $f_{\Lambda_0}(p)$
of the form $e^{\alpha/2}(w)$ for some $w \in W(\Lambda_1)$? 
This is the case if and only if
for some $q=q(y_{-2},y_{-3},\dots) \in \mathcal{A}'$,
\bea \label{10.17}
&& f_{\Lambda_0}(p)= e^{\alpha/2} (f'_{\Lambda_1}(q))= e^{\alpha/2}
(f_{\Lambda_1}(q))=\nn
&& e^{\alpha/2}(q(x_{\alpha}(-2),x_{\alpha}(-3),\ldots )
e^{\alpha/2}v_{\Lambda_0}=
q(x_{\alpha}(-3),x_{\alpha}(-4),\ldots ) e^{\alpha}v_{\Lambda_0}=\nn
&& q(x_{\alpha}(-3),x_{\alpha}(-4),\ldots ) x_{\alpha}(-1)v_{\Lambda_0}
= f_{\Lambda_0}(q(y_{-3},y_{-4},\ldots ) y_{-1})= \nn
&&= f_{\Lambda_0}((\widehat {e^{\alpha/2}}(q))y_{-1}).
\eea
This in turn is the case if and only if for some $q\in \mathcal A'$
$$
p-q(y_{-3},y_{-4},\ldots ) y_{-1}\in \mathcal A_{\Lambda_0},
$$
i.e., if and only if 
$$
p\in \widehat {e^{\alpha/2}}( \mathcal A')y_{-1}+ \mathcal
A_{\Lambda_0}.
$$

Note that
$$
\widehat {e^{\alpha/2}}( \mathcal A')y_{-1} = 
\mathbb{C}[y_{-3},y_{-4},\ldots ]y_{-1}.
$$
Thus 
$$
f_{\Lambda_0}(p)\in {\rm Im} \ e^{\alpha/2}: \ W(\Lambda_1) \rightarrow W(\Lambda_0)
 \iff p\in \mathbb{C}[y_{-3},
y_{-4},\ldots ]y_{-1}+ \mathcal A_{\Lambda_0}.
$$

Now we prove that in fact
\be \label{10.22}
f_{\Lambda_0}(p)\in {\rm Im} \ e^{\alpha/2}: \ W(\Lambda_1) \rightarrow W(\Lambda_0)
\iff p\in (y_{-1})+ \mathcal
A_{\Lambda_0}.
\ee
For this, what we need to prove is that
\be \label{10.23}
(y_{-1})\subset \mathbb{C}[y_{-3},y_{-4},\ldots ]y_{-1}+
\mathcal A_{\Lambda_0},
\ee
and for this we use Proposition \ref{x20}:  We have
$$
X(\alpha,x)^2=0
$$ 
on $L(\Lambda_0)$ (and on $L(\Lambda_1)$), i.e., for all $j\in \mathbb{Z}$,
$$
\sum_{i\in \mathbb{Z}} x_{\alpha}(j-i)x_{\alpha}(i)=0.
$$
In particular, applying these relations to
$v_{\Lambda_0}$, we see that 
for all $n\leq -2$,
\be \label{10.26}
\sum_{i,j<0,\ i+j=n} x_{\alpha}(i)x_{\alpha}(j)v_{\Lambda_0}=0.
\ee
The first two such relations are
$$
x_{\alpha}(-1)^2v_{\Lambda_0}=0
$$
and
$$
x_{\alpha}(-2)x_{\alpha}(-1)v_{\Lambda_0}=0.
$$
That is,
$$
y_{-1}^2\in \mathcal A_{\Lambda_0}
$$
and
$$
y_{-2}y_{-1}\in \mathcal A_{\Lambda_0}.
$$
Thus any monomial in $\mathcal A$ divisible by $y_{-1}^2$ or by
$y_{-2}y_{-1}$ lies in
$\mathcal A_{\Lambda_0}$, and this proves (\ref{10.23}) and hence
(\ref{10.22}). 
Thus for $p \in \mathcal{A}$,
\bea \label{10.31}
&& f_{\Lambda_0}(p)\in {\rm Im}\ e^{\alpha/2} \iff p\in (y_{-1})+ \mathcal
A_{\Lambda_0} \nn
&& \iff \pi(p) \in \pi(\mathcal A_{\Lambda_0}) \iff S(p)\in 
S(\mathcal A_{\Lambda_0})
\eea
(where we use the definition (\ref{10.11}) of $S$), 
and we had already seen that
\bea \label{10.32}
&& f_{\Lambda_0}(p)\in {\rm Ker}\ o(e^{\alpha/2})\iff p\in (y_{-1}) \oplus \mathcal
B_{\Lambda_1} \nn
&& \iff \pi(p)\in \mathcal B_{\Lambda_1} \iff S(p)\in 
\mathcal A_{\Lambda_0}.
\eea

We shall finally use these characterizations of ${\rm Ker} \
o(e^{\alpha/2})$ and ${\rm Im} \ e^{\alpha/2}$ to prove the exactness.

First we confirm the chain-complex property, which we had already seen:
$$
{\rm Im} \ e^{\alpha/2}\subset {\rm Ker} \ o(e^{\alpha/2})
$$
since
\be \label{10.34}
\pi(\mathcal A_{\Lambda_0})\subset \mathcal B_{\Lambda_1}
\ee
(recall (\ref{10.9})) or since
$$S(\mathcal{A}_{\Lambda_0}) \subset \mathcal{A}_{\Lambda_0}$$
(recall (\ref{10.16})).

Now we prove exactness. 
What we must prove is that 
$$
\mathcal B_{\Lambda_1} \subset \pi(\mathcal A_{\Lambda_0})
$$
or that
$$
\mathcal A_{\Lambda_0}\subset S(\mathcal A_{\Lambda_0}),
$$
i.e., that
$$
\pi(\mathcal A_{\Lambda_0})= \mathcal B_{\Lambda_1}
$$
or that
\be \label{10.38}
S(\mathcal A_{\Lambda_0})=\mathcal A_{\Lambda_0}.
\ee

\begin{remark} 
This is a precise statement of the principle that ``the only difference
in the relations defining $W(\Lambda_0)$
and $W(\Lambda_1)$ is the initial condition 
relation $x_{\alpha}(-1)v_{\Lambda_1}=0$''.
But it is expressed as relations involving only
$\mathcal{A}_{\Lambda_0}$ and not $\mathcal{A}_{\Lambda_1}$
or $\mathcal{B}_{\Lambda_1}$:
$$
S(\mathcal A_{\Lambda_0})=\mathcal A_{\Lambda_0}.
$$
\end{remark}

To prove this,
consider the relations (\ref{10.26})
in $W(\Lambda_0)$.
These relations amount to the assertion that
$$
r_n= \sum_{i,j<0, \ i+j=n} y_{i}y_{j}\in \mathcal A_{\Lambda_0}
$$
for $n\leq -2$. Note that $r_{-2}=y_{-1}^2$ and
$r_{-3}=y_{-2}y_{-1}$, which we discussed above.
By Feigin-Stoyanovsky \cite{FS1}, these 
elements $r_n$ {\it generate}
the ideal $\mathcal A_{\Lambda_0}$ in $\mathcal A$, that is,
\be \label{10.42}
\mathcal A_{\Lambda_0}= \sum_{n \leq -2}\mathcal A\ r_n.
\ee
{}From this, it is clear that $S(\mathcal A_{\Lambda_0})=\mathcal A_{\Lambda_0}$, since
$$
S(r_n)=r_{n+2}
$$
for $n\leq -4$ and 
$$
S(r_{-3})=S(r_{-2})=0.
$$
This proves the exactness.

The main theorem immediately gives:

\begin{theorem} 
We have the relation
\be \label{10.45}
\chi_0(x,q)= x^{-1/2}q^{-1/4} \chi_1(x,q)+ x^{1/2}q^{1/4}{\chi}_1(xq,q).
\ee
\end{theorem}
{\em Proof:} The operator $e^{\alpha/2}$ has weight equal to the charge of
the source vector plus $1/4$,
and it increases charge by $1/2$. On the other hand $o(e^{\alpha/2})$
is of weight
$1/4$ and it increases the charge by $1/2$. 
The exactness of (\ref{10.1}) now immediately gives (\ref{10.45}).

We now have the Rogers-Ramanujan recursion:

\begin{corollary} 
$$
\chi_0(x,q)= \chi_0(xq,q)+ xq\chi_0(xq^2,q).
$$
\end{corollary}
{\em Proof:}
Apply formula (\ref{9.36}) and the Theorem.

Solving this recursion (cf. \cite{A}), we immediately obtain:

\begin{corollary}
\bea
&& \chi_{0}(x,q)=\sum_{n\geq 0}
\frac{x^nq^{n^2}}{(q)_n},\nn
&& \chi_{1}(x,q)= x^{1/2}q^{1/4}\sum_{n\geq 0}
\frac{x^n q^{n^2+n}}{(q)_n}. \nonumber
\eea
\end{corollary}

\noindent {\small \sc 
Dipartimento Me. Mo. Mat.,
Universit{\`a} di Roma ``La Sapienza,''
Via A. Scarpa 16,
00161 Roma, Italia} \\
{\em E-mail address}: capparelli@dmmm.uniroma1.it \\
\noindent {\small \sc Department of Mathematics,  Rutgers University, Piscataway,
NJ 08854} \\
{\em E-mail address}: lepowsky@math.rutgers.edu \\
{\small \sc Department of Mathematics, University of Arizona, Tucson,
AZ 85721} \\
{\em E-mail address}: milas@math.arizona.edu

\end{document}